\documentstyle[amscd,amssymb]{article} 

\makeatletter
\@ifundefined {upshape}{\newcommand {\upshape}{\fontshape {n}\selectfont}}{}
\makeatother

\newcommand{\im}{\mathop{\mbox{\upshape Im}}\nolimits}
\newcommand{\Aut}{\mathop{\mbox{\upshape Aut}}\nolimits}
\newcommand{\Alg}{\mathop{\mbox{\upshape Alg}}\nolimits}
\newcommand{\Hdg}{\mathop{\mbox{\upshape Hdg}}\nolimits}

\newcommand{\Sing}{\mathop{\mbox{\upshape Sing}}\nolimits}
\newcommand{\lspan}{\mathop{\mbox{\upshape span}}\nolimits}
\newcommand{\into}{\hookrightarrow}
\newcommand{\PR}{{\Bbb P}}
\newcommand{\cO}{{\cal O}}
\newcommand{\cZ}{{\cal Z}}
\newcommand{\C}{{\Bbb C}}
\newcommand{\Z}{{\Bbb Z}}
\newcommand{\dsdxj}{\frac{\partial s}{\partial x_j}}
\newcommand{\rk}{\mathop{\mbox{\upshape rk}}\nolimits}
\newcommand{\implies}{\Longrightarrow}

\newenvironment{pf}{\noindent {\sc Proof:\/}}{\hfill $\boxtimes$}
\newtheorem{thm}{Theorem}

\newtheorem{pro}[thm]{Proposition}
\newtheorem{lem}[thm]{Lemma}
\newtheorem{rmk}[thm]{Remark}

\begin{document}
\pagenumbering{roman}
\title{A Noether-Lefschetz theorem for vector bundles}
\author{Jeroen G.~Spandaw}
\date{July 7, 1995}
\maketitle

\begin{abstract} In this note we use the monodromy argument to
prove a Noether-Lefschetz theorem for vector bundles.
\end{abstract}

\section{Introduction}\pagenumbering{arabic}

Let $X$ be a smooth complex projective manifold of dimension $n$ and let $E$
be a very ample vector bundle on $X$ of rank $r$. This means that the
tautological quotient line bundle $L$ on the bundle
$Y=\PR(E^\ast)$ of hyperplanes in $E$ is very ample.
For almost all $s\in H^0(X,E)$ the zero-locus
$Z$ is smooth, irreducible and of dimension $n-r$.
In \cite[prop.~1.16]{S} Sommese
proved that $H^{i}(X,Z;\Z)$ vanishes for $i<n-r+1$ and is torsion free for
$i=n-r+1$. Assume that $n-r$ is even, say $n-r=2p$.
Let $\Alg\subset H^{n-r}(Z)$ be the space of algebraic classes
and let $\im=\im(H^{n-r}(X)\into H^{n-r}(Z))$.
(We always take coefficients in $\C$
unless other coefficients are mentioned explicitely (cf.~Remark~\ref{rmk}).)
In this note we prove the following Noether-Lefschetz theorem for this
situation.

\begin{thm}\label{thm1} If $E$ is very ample and $s$ is general,
then either $\Alg\subset \im$ or $\Alg+\im=H^{n-r}(Z)$.
\end{thm}
(With \lq\lq general\rq\rq we shall always mean general in the usual
Noether-Lefschetz sense.)
The following theorem, which generalizes the Noether-Lefschetz
theorem for complete intersections in projective space
(see \cite[pp.~328--329]{DK})
is an immediate corollary.

\begin{thm}\label{thm2} If $h^{\alpha\beta}(X)<h^{\alpha\beta}(Z)$ for some
pair $(\alpha,\beta)$ with
$\alpha+\beta=n-r$ and $\alpha\neq \beta$, then every algebraic class on
$Z$ is induced
from $X$.
\end{thm}

\begin{rmk}\label{rmk}\normalshape
Notice that the unique pre-images of
algebraic classes are themselves Hodge classes, i.e.\ lie in
$H^{p,p}(X)\cap H^{n-r}(X;\Z)$. This follows from the fact that the
cokernel  of $H^{n-r}(X,\Z)\to H^{n-r}(Z,\Z)$ is torsion free.
\end{rmk}

It is not difficult to show that
after replacing $E$ with $E\otimes L^k$, where
$k\gg0$ and $L$ is an ample line bundle,
the assumption of theorem~\ref{thm2} is satisfied.
(E.g.\ the geometric genus of $X$ goes to infinity as
$k$ goes to infinity.)
In \cite{Sp} we used the notion of Castelnuovo-Mumford regularity
(cf.~\cite[p.~99]{M}) to make
the positivity assumption on $E$ more precise
if $X=\PR^n$.
Notations are as in theorem~\ref{thm1}.
$\Hdg$ is defined to be the space of Hodge classes
on $Z$ of codimension~$p$,
i.e.\ $\Hdg=H^{p,p}(Z)\cap H^{n-r}(Z,\Z)$.

\begin{thm}\label{thm3} If $E$ is a $(-3)$-regular vector bundle of rank
$r$ on $X=\PR^n$ and
$Z$ is the zero-locus of a general global section of $E$, then
$\Hdg\subset\im$,
unless $(X,E)=(\PR^3,\cO(3))$.
If $\dim Z=2$, then it suffices that $E$ be $(-2)$-regular
unless $(X,E)=(\PR^3,\cO(2))$, $(\PR^3,\cO(3))$ or
$(\PR^4,\cO(2)\oplus\cO(2))$.
\end{thm}
(Notice that $(-3)$-regularity $\implies$ $(-1)$-regularity
$\implies$ very ampleness.)
For the case $\dim Z=2$ theorem~\ref{thm3} is due to
Ein \cite[thm.~3.3]{E}.
The advantage of theorem~\ref{thm3}
is that it applies to {\em Hodge\/} rather than {\em
algebraic\/} classes on $Z$. For example, it implies that
if all Hodge classes of codimension $n-r$ on $\PR^n$ are
algebraic, then the same holds for $Z$.
The advantage of theorem~\ref{thm1} is that the
positivity condition on $E$ is more geometric:
the cohomological conditions from
\cite{Sp} are replaced with the condition that $E$ be very ample plus
a Hodge number inequality (cf.~theorem~\ref{thm2}).
In other words, for very ample vector bundles, the Noether-Lefschetz
property holds as soon as this is allowed by the Hodge numbers.
However, this Hodge number inequality condition
is of course a cohomological condition on $E$ in disguise.

\smallskip

\noindent {\em Acknowledgement\/} I am grateful to professor Sommese for
the suggestion
that I look at the bundle $\pi\colon \PR(E^\ast)\to X$
of hyperplanes in $E$.

\section{Proof of the main result}

Let $V=H^0(X,E)$,
let $\PR(V)$ be the set of lines in $V$,
let $N=\dim\PR(V)=h^0(X,E)-1$
and set $X'=\PR(V)\times X$.
Set $E'=p_1^\ast\cO(1)\otimes p_2^\ast E$,
where $p_i$ are the projections.
$E'$ has a canonical section $s'$.
Let $\cZ$ be the zero locus of $s'$.
The restriction
$
     p\colon\cZ\to\PR(V)
$
of $p_1$ to $\cZ$ is
the universal family of
zero loci of sections in $E$.
We leave the proof of the following easy lemma to the reader.
\begin{lem} If $E$ is very ample, then it is generated by its sections.
If $E$ is generated by its sections, then $\cZ$ is smooth, irreducible
and of dimension $N+n-r$.
\end{lem}

Let $\Delta\subset\PR(V)$ be the discriminant of $p$, i.r.\
\begin{eqnarray*}
     \Delta&=&p\{z\in\cZ: \rk_z p\le N-1\}\\
&=&\{[s]\in\PR(V): \hbox{$p^{-1}(s)$ is not smooth of dimension $n-r$}\}.
\end{eqnarray*}
Fix a point $[s_0]\in\PR(V)\setminus\Delta$ and let $Z\subset X$
be the corresponding smooth fibre of $p$.
Let $\Gamma$ the image of the monodromy representation
$\pi_1(\PR(V) \setminus\Delta)\to \Aut(H^{n-r}(Z))$.

Let $\im^\perp$ be the orthogonal complement of $\im$ with
respect to the intersection form on $H^{n-r}(Z)$.
Since for general $s\in H^0(X,E)$, $\Alg$ is a $\Gamma$-module
(cf.\ \cite[p.~141]{H}),
theorem~\ref{thm1} from the following proposition.

\begin{pro} (Second Lefschetz Theorem)
\begin{enumerate}
\item $H^{n-r}(Z)=\im\oplus\im^\perp$
\item $\im=H^{n-r}(Z)^\Gamma$
\item $\im^\perp$ is an irreducible $\Gamma$-module
\end{enumerate}
\end{pro}
\begin{pf}
\begin{enumerate}
\item Arguing as in the proof of \cite[thm.~6.1 (i)]{G}
one shows that if $Z$ is submanifold of
a compact K\"ahler manifold $X$ such that
$H^{i}(X,Z)=0$ for $i\le m= \dim Z$,
then the restriction of the intersection form
to $\im(H^m(X)\into H^m(Z))$ is non-degenerate.
\item The inclusion $\im\subset H^{n-r}(Z)^\Gamma$ is trivial. To prove that
$H^{n-r}(Z)^\Gamma\subset\im$,
we argue as in \cite[thm.~6.1 (iii)]{G}.
Consider the commutative diagram
$$
\begin{CD}
H^{n-r}(\PR(V)\times X) @>>> H^{n-r}(\cZ)\\
@VVV @VVV\\
H^{n-r}(X) @>>> H^{n-r}(Z)^\Gamma.
\end{CD}
$$
By \cite[th\'eor\`eme 4.1.1 (ii)]{D} the map
$H^{n-r}(\cZ)\to H^{n-r}(Z)^\Gamma$ is surjective.
By \cite[prop.~1.16]{S} the map
$H^{n-r}(\PR(V)\times X) \to H^{n-r}(\cZ)$
is surjective.
\item Since the monodromy respects the intersection form,
$I^\perp$ is a $\Gamma$-module.
The standard argument using Lefschetz pencils
and the theory of vanishing cycles reduces
the problem of irreducibility to proposition~\ref{prop} below
(cf. \cite[pp.~46--48]{L}).
\end{enumerate}
\end{pf}

\begin{pro}\label{prop}
\begin{enumerate}
\item The discriminant $\Delta$ is an irreducible,
closed, proper subvariety of $\PR(V)$.
\item Let $G\subset\PR(V)$ be a general line.
Then $\cZ_G:=p^{-1}(G)$ is smooth, irreducible of dimension $n-r+1$
and the restricted family $p_G\colon\cZ_G\to G$ is a holomorphic
Morse function, i.e. all critical points are non-degenerate
and no two lie in the same fibre (cf.~\cite[p.~34]{L}).
$g\in G$ is a critical value of $p_G$ if and only if it
is a critical value of $p$.
\end{enumerate}
\end{pro}
\begin{pf} The statements about $\cZ_G$ follow
from Bertini. The remaining assertions
are well-known if $\rk E=1$ (cf.~\cite[p.~19]{L}).
In particular, they are true
for $(Y,L)$, where $Y$ is the hyperplane bundle $\PR(E)$ of $E$
and $L$ is the tautological quotient line bundle $\cO_Y(1)$.
The following proposition reduces the general case $(X,E)$
to this line bundle case $(Y,L)$, thus finishing the proof.
\end{pf}

Before we state the last proposition, notice that
the natural map $s\mapsto \bar{s}\colon H^0(X,E)\to H^0(Y,L)$,
where $\bar{s}(x,h):=\overline{s(x)}\in E(x)/h=L(x,h)$
for $(x,h)\in Y$, is an isomorphism. Indeed, the map
is clearly injective and $h^0(Y,L)=h^0(X,\pi_\ast L)=h^0(X,E)$.
For $s\in H^0(X,E)$ we denote by $Z_X(s)$
the zero-locus of $s$ in $X$
and by $Z_Y(\bar{s})$ the zero-locus of $\bar{s}$ in $Y$.
\begin{pro} For $s\in H^0(X,E)\setminus\{0\}$,
$Z=Z_X(s)$ is singular
if and only if $W=Z_Y(\bar{s})$ is singular.
More precisely, if $x\in\Sing Z$, then there exists a $y\in\Sing W$
with $\pi(y)=x$ and conversely,
if $(x,h)\in\Sing Z$, then $x\in \Sing W$.
Finally, if $(x,h)$ is a non-degenerate quadratic singularity,
then so is $x$.
\end{pro}
\begin{pf} This is a calculation in local coordinates.
Let $x_0\in Z$, i.e.\ $s(x_0)=0$.
After choosing local coordinates $x_1,\ldots,x_n$ on $X$
and a local trivialization of $E$ near $x_0$
we may regard $s$ to be a function in $x_1,\ldots,x_n$.
Then $x_0\in\Sing Z$ if and only if
$\{\dsdxj(x_0)\}_{j=1}^n$ does not span $\C^r$.
Let $h_0\subset\C^r$ be a hyperplane containing
$\lspan\{\dsdxj(x_0)\}_{j=1}^n$. We claim that $y_0=(x_0,h_0)\in\Sing W$.
We may assume that the local trivialization of $E$
has been chosen in such a way that $h_0$ is given
by $z_r=0$, where $z_1,\ldots,z_r$ are coordinates on $\C^r$.
Let $s=(f_1,\ldots,f_r)$.
Local coordinates on $Y$ near $y_0$ are provided by
the local coordinates $x_1,\ldots,x_n$ on $X$ near $x_0$
together with $(y_1,\ldots,y_{r-1})\in\C^{r-1}$:
we let $(y_1,\ldots,y_{r-1})\in\C^{r-1}$ correspond to
the hyperplane $\sum_{i=1}^r y_iz_i=0$,
where $y_r:=1$. The point $y_0$ has coordinates $(x_0,0)$.
In these local coordinates $\bar{s}(x,y)=\sum_{i=1}^ry_if_i(x)$.
It now suffices to calculate $\frac{\partial \bar{s}}{\partial x_k}(x_0,0)
=\frac{\partial f_r}{\partial x_k}(x_0)=0$
for $k=1,\ldots,n$
and $\frac{\partial \bar{s}}{\partial y_j}(x_0,0)=f_j(x_0)=0$
for $j=1,\ldots,r-1$.
The converse is proven similarly.

Let $y_0=(x_0,h_0)\in\Sing W$. We may again assume that
$h_0$ is given by $z_r=0$.
The Hessian of $\bar{s}$ in $y_0$ is of the form
$\left(\begin{array}{cc} h & d^t\\ d & 0\end{array}\right)$,
where the $n\times n$-matrix $h$ is the Hessian of $f_r$
and the $(r-1)\times n$-matrix $d$ is the Jacobian
of $f':=(f_1,\ldots,f_{r-1})$ in $x_0$.
Let $Z'=\{x\in X: f'(x)=0\}$.
We have to check that the Hessian of $f_r|_{Z'}$ in $0$ is non-degenerate.
Since we assume that
the Hessian of $\bar{s}$ has maximal rank in $y_0$, so has $d$. Thus,
after a change of coordinates, we may assume that $f_i(x)=x_i$ for $i<r$.
Then $\bar{s}(x,y)=\sum_{i=1}^{r-1}x_iy_i+f_r(x)$,
hence the Hessian of $\bar{s}$ in $y_0$ is
$$
\left(\begin{array}{ccc}
\ast & \ast & E_{r-1}\\
\ast & H & 0\\
E_{r-1} & 0 & 0,
\end{array}
\right),
$$
where $H$ is the Hessian of $f_r|_{Z'}$ in $x_0$.
It follows that $H$ is non-degenerate.
\end{pf}

\noindent Jeroen Spandaw\\
Institut f\"ur Mathematik\\
Universit\"at Hannover\\
Postfach 6009\\
D-30060 Hannover\\
Germany\\
e-mail: spandaw@math.uni-hannover.de

\end{document}